\documentclass[12pt]{amsart}
\usepackage{amsmath}
\usepackage{amsfonts}
\usepackage{amssymb}
\usepackage{amscd}

\newcommand{\be}{\beta}
\newcommand{\ga}{\gamma}
\newcommand{\de}{\delta}

\newcommand{\e}{\varepsilon}

\newcommand{\Si}{\Sigma}
\newcommand{\BR}{\mathbb{R}}
\newcommand{\BZ}{\mathbb{Z}}
\newcommand{\BQ}{\mathbb{Q}}
\newcommand{\BN}{\mathbb{N}}
\newcommand{\A}{\mathcal{A}}

\newcommand{\wt}{\widetilde}

\newcommand{\per}{\mathrm{Per}}

\newcommand{\n}{\mathfrak{n}}

\newcommand{\BT}{\mathbb{T}}

\renewcommand{\A}{\mathcal{A}}
\renewcommand{\phi}{\varphi}

\newtheorem{lemma}{Lemma}

\newtheorem{thm}[lemma]{Theorem}
\newtheorem{cor}[lemma]{Corollary}
\theoremstyle{definition}
\newtheorem{Def}[lemma]{Definition}

\theoremstyle{remark}
\newtheorem{rmk}[lemma]{Remark}

\begin{document}

\title[Ergodic-theoretic properties of Bernoulli
convolutions] {Ergodic-theoretic properties of \\ certain
Bernoulli convolutions}
\author{Nikita Sidorov}
\address{Department of Mathematics, UMIST, P.O. Box 88,
Manchester M60 1QD, United Kingdom. E-mail:
Nikita.A.Sidorov@umist.ac.uk}
\date{\today}
\thanks{Supported by the EPSRC grant no GR/R61451/01.}
\subjclass[2000]{28D05, 11R06} \keywords{Bernoulli convolution,
$\be$-shift, Erd\"os measure, ergodic properties}
\begin{abstract}
In \cite{SV} the author and A.~Vershik have shown that for
$\be=\frac12(1+\sqrt5)$ and the alphabet $\{0,1\}$ the infinite
Bernoulli convolution ($=$ the Erd\"os measure) has a property
similar to the Lebesgue measure. Namely, it is quasi-invariant
of type $\mathrm{II}_1$ under the $\be$-shift, and the natural
extension of the $\be$-shift provided with the measure
equivalent to the Erd\"os measure, is Bernoulli. In this note we
extend this result to all Pisot parameters $\be$ (modulo some
general arithmetic conjecture) and an arbitrary ``sufficient"
alphabet.
\end{abstract}

\maketitle

\noindent\textbf{1. Introduction and the main theorem.} Let
$\be>1$; the \textit{infinite Bernoulli convolution} (or the
\textit{infinitely convolved Bernoulli measure}) is defined as the
infinite convolution of the independent discrete random variables
$\theta_n(\be)$ for $n$ from 1 to $\infty$, where $\theta_n(\be)$
assumes the values $\pm\be^{-n}$ with the probability $\frac12$.
This measure is well studied from the probabilistic point of view
-- see, e.g., \cite{Ga, AlZa}. In particular, if $\be>2$, then the
support of the corresponding infinite Bernoulli convolution is a
Cantor set of zero Lebesgue measure, and for $\be=2$ it coincides
with the Lebesgue measure on $[-1,1]$. Besides, if $\be$ is a
\textit{Pisot number} (i.e., an algebraic integer $>1$ whose
Galois conjugates are all less than 1 in modulus), then the famous
Erd\"os Theorem claims that it is singular with respect to the
Lebesgue measure \cite{E}. Finally, it is also worth mentioning
the fundamental result by B.~Solomyak who has proved that it is
absolutely continuous for a.e. $\be\in(1,2)$ \cite{So}. The aim of
this short note is to study some ergodic-theoretic properties of
this important measure in the case of Pisot parameter $\be$.

Actually, we will consider a slightly more general model. Namely,
let $d\in\BN\setminus\{1\}, \A_d=\{0,1,\dots,d-1\}$ and $\be$ be
an irrational Pisot number, $1<\be<d$.

Let $\mu^+_{\be,d}$ denote the {\it Erd\"os measure}, i.e., the
measure on $\BR$ that corresponds to the distribution of the
random variable
\begin{equation}
\xi_{\be,d}=\pi_{\be,d}((x_n)_1^\infty):=
\frac{\be-1}{d-1}\sum_{n=1}^\infty x_n\be^{-n}, \label{pr}
\end{equation}
where $x_n$'s are i.i.d. variables, each of which assumes the
values $\{0,1,\linebreak[0]\dots,d-1\}$ with the probability
$1/d$. Since $\be<d$, it is obvious that
$\textrm{supp}\,\mu_{\be,d}=[0,1]$. Let $\tau^+_\be$ denote the
$\be$-shift in $[0,1)$, i.e.,
$$
\tau^+_\be(x)=\be x\bmod1.
$$
The relationship between the Erd\"os measure and the infinite
Bernoulli convolutions is straightforward: let $d=2$ and
$\be\in(1,2)$. Then the affine map
$h_\be(x)=\frac{\be-1}2x+\frac12$ turns the corresponding infinite
Bernoulli convolution into the Erd\"os measure with the same
parameter $\be$. Since an affine transform does not alter any
essential ergodic properties, we may confine ourselves to the
study of the measures $\mu^+_{\be,d}$.

Let $X_\be^+$ denote the one-sided {\em $\be$-compactum}, i.e.,
the space of all possible (greedy) $\be$-expansions of the numbers
in $[0,1)$.

More precisely, let the sequence $(a_n)_1^\infty$ be defined as
follows: let $1=\sum_{1}^{\infty}a_k' \be^{-k}$ be the greedy
expansion of 1, i.e, $a_n'=[\be(\tau^+_\be)^{n-1}1],\ n\ge1$. If
the tail of the sequence $(a_n')$ differs from $0^\infty$, then we
put $a_n\equiv a_n'$. Otherwise let $k=\max\,\{j:a_j'>0\}$, and
$(a_1,a_2,\dots):=
(a_1',\dots,a_{k-1}',\linebreak[0]a_k'-1)^\infty$. In the seminal
paper \cite{Pa} it is shown that for each greedy expansion $\e$ in
base $\be,\ (\e_n,\e_{n+1},\dots)$ is lexicographically less
(notation: $\prec$) than $(a_1,a_2,\dots)$ for every $n\ge1$.
Moreover, it was shown that, conversely, every sequence with this
property is actually the greedy expansion in base $\be$ for some
$x\in [0,1)$.

Put
\[
X_\be^+=\left\{\e\in\prod_1^\infty\{0,1,\dots,[\be]\}\mid
(\e_n,\e_{n+1},\dots)\prec(a_1,a_2,\dots),\ n\in\BN\right\}
\]
(the one-sided $\be$-compactum), and
\[
X_\be=\left\{\e\in\prod_1^\infty\{0,1,\dots,[\be]\}\mid
(\e_n,\e_{n+1},\dots)\prec(a_1,a_2,\dots),\ n\in\BZ\right\}
\]
(the two-sided $\be$-compactum). The sequences from the
$\be$-compactum (one-sided or two-sided) will be called {\em
$\be$-expansions}. It follows from the above formulas that both
$\be$-compacta are stationary ($=$ shift-invariant). As was shown
in \cite{Pa}, the map $\phi_\be:X_\be^+\to[0,1)$ defined by the
formula
\begin{equation}
\phi_\be(\e)=\sum_{n=1}^\infty\e_n\be^{-n}\label{bexp},
\end{equation}
is one-to-one except a countable set of sequences.

Let $Fin(\be)$ denote all numbers from $[0,1)$ that have a finite
$\be$-expansion (i.e., the tail is $0^\infty$). It is obvious that
$Fin(\be)\subset\BZ[\be]\cap[0,1)$; the inverse inclusion however
is not necessarily true even for a Pisot number, see, e.g.,
\cite{Ak, S1}.

\begin{Def}\label{WF} We call a Pisot number $\be$ {\em weakly
finitary} (WF) if for any $y\in\BZ[\be]\cap[0,1)$ and any
$\delta>0$ there exists $f\in Fin(\be)\cap(0,\delta)$ such that
$y+f\in Fin(\be)$ as well.
\end{Def}

\begin{rmk}
The notion of WF number appeared in different settings in a number
of recent works \cite{Ak, S1, S2} and in a slightly different form
-- earlier in the thesis \cite{Hol}. There exists a conjecture
(shared by most experts in the area) that in fact {\bf every}
Pisot number is weakly finitary. Note that Sh.~Akiyama \cite{Ak}
has given an explicit algorithm of checking whether a {\em given}
Pisot number is WF, and, as far as we are concerned, none of them
has failed so far.
\end{rmk}

The main theorem of the present note is as follows.

\begin{thm}\label{main} If $\be$ is WF, then the Erd\"os measure
$\mu^+_{\be,d}$ is quasi-invariant under the $\be$-shift
$\tau^+_\be$. Moreover, there exists a unique probability measure
$\nu^+_{\be,d}$ invariant under $\tau^+_\be$ and equivalent to
$\mu^+_{\be,d}$. The natural extension of the endomorphism
$([0,1), \nu^+_{\be,d},\tau^+_\be)$ is Bernoulli.
\end{thm}

We believe that knowing this fact could be important for the study
of further ergodic-theoretic properties of this important measure,
including its Gibbs structure and multifractal spectrum (see
\cite{OST} for some results in this direction and references
therein).

\medskip\noindent\textbf{2. Auxiliary results and definitions.}
The rest of the paper is devoted to the proof of
Theorem~\ref{main}, which is based on the idea of \cite[\S1]{SV}
(where the special case $\beta=\frac{1+\sqrt5}2,\ d=2$ was
considered) and also uses techniques of \cite{S2}. Note that the
claim analogous to Theorem~\ref{main} is known to be true for the
Lebesgue measure on $[0,1]$ -- see \cite{Sm, DKS}. The above
theorem therefore immediately leads to an ergodic-theoretic proof
of the famous Erd\"os Theorem which claims that the Erd\"os
measure is singular \cite{E}. Indeed, it suffices to apply the
corollary of the Birkhoff Ergodic Theorem claiming that two
ergodic measures either coincide or are mutually singular; the
fact that the two invariant measures in question do not coincide
can be proved in the very same way as in the case
$\be=\frac{1+\sqrt5}2$, see \cite[Proposition~1.10]{SV}.

Our first goal is to define the two-sided normalization in base
$(\be,d)$. Let $\Si_d:=\prod_{-\infty}^\infty\{0,1,\dots,d-1\},
\Si_d^+:=\prod_1^\infty\{0,1,\dots,d-1\}$. We will use the
following convention: the sequences from $X_\be$ ($X_\be^+$) will
be denoted with the letter ``$\e$" and sequences from the full
compacta -- with ``$x$". Let $p_d$ denote the product measure on
$\Si_d$ with the equal multipliers and $p_d^+$ -- its one-sided
analog. Recall that the {\em one-sided normalization}
$\n_{\be,d}^+$ is defined as the map from $\Si_d^+$ to $X_\be^+$
acting by the formula
\begin{equation}
\n_{\be,d}^+=\phi_\be^{-1}\circ\pi_{\be,d}, \label{n+}
\end{equation}
where $\pi_{\be,d}$ is given by (\ref{pr}) and $\phi_\be$ is given
by (\ref{bexp}) -- see \cite{Fr}. The following convention will be
used hereinafter: the notation $\n_{\be,d}^+(x_1\dots x_n)$ means
$\n_{\be,d}^+(x_1,\dots,x_n,0^\infty)$ and if
$\n_{\be,d}^+(x_1\dots x_n)=(\e_1,\dots,\e_{n'},0^\infty)$, then
by definition, $\n_{\be,d}^+(x_1\dots x_n)=(\e_1,\dots,\e_{n'})$,
i.e., we ignore the tail $0^\infty$ whenever possible. By the
above, the Erd\"os measure may be computed by the formula
$$
\mu^+_{\be,d}=\n^+_{\be,d}(p_d^+).
$$
For more details see \cite{SV}.

\begin{rmk} There exists a more direct way of defining
normalization. Namely, in \cite{Fr} it was shown that one may find
a finite automaton that carries out the operation of normalization
in Pisot bases. The converse is also true: if the function of
normalization is computable by a finite automaton, then $\be$ must
be a Pisot number -- see \cite{BeFr}.
\end{rmk}


We need one more technical lemma  before we may proceed.

\begin{lemma}\label{L}
If $\be$ is WF, then there exists $L=L(\be,d)\in\BN$ such that
for any word $x_1\dots x_n$ in the alphabet $\A_d$ there exists
a word $x_{n+1}\dots x_{n+L}$ in the same alphabet such that
$\n^+_{\be,d}(x_1\dots x_{n+L})$ is finite.
\end{lemma}
\begin{proof} By the well-known result of K.~Schmidt \cite{Sch},
the $\be$-expansion of any $x\in\BQ[\be]\cap(0,1)$ is eventually
periodic; moreover, for the elements of $\BZ[\be]\cap(0,1)$ the
collection of such periods is known to be finite \cite{Ak, S2}.
Let $\mathcal{T}_\be=(\per_1,\dots,\per_r)$ denote this
collection. In \cite{FrSo} it was shown that the normalization of
a word $x_1\dots x_n$ in base $(\be,d)$ has the following form:
\[
\n_{\be,d}^+(x_1\dots
x_n)=(\e_1,\dots,\e_n,\e_{n+1},\dots,\e_{n+L_1},\per_j),\quad 1\le
j\le r,
\]
where $\per_j\in\mathcal T_\be$ and $L_1$ is a function of $\be$
and $d$. Thus, it suffices to prove the claim for the words $x$ of
the form $\e_1\dots\e_{n+L_1}x^{(j)}_1\dots x^{(j)}_{p_j}$, where
$x^{(j)}_1\dots x^{(j)}_{p_j}$ is the word in the alphabet $\A_d$
whose normalization is $\per_j$. Let $p=\max_{j=1}^r p_j,\
\de=\frac{d-1}{\be-1}\,\be^{-n-L_1-p}$ and
$y=\sum_{i=1}^{n+L_1}\e_i\be^{-i}+
\sum_{i=1}^{p_j}x_i^{(j)}\be^{-i-n-L_1}$. Then by
Definition~\ref{WF}, there exists $f\in Fin(\be)\cap(0,\de)$ such
that $y+f\in Fin(\be)$. By our choice of $\de$, the
$\be$-expansion of $f$ must be of the form $(f_{n+L_1+p+1},\dots,
f_{n+L_1+p+L_2})$ for some fixed $L_2$ (because we have only a
finite number of $(x^{(j)}_1\dots x^{(j)}_{p_j})$). Setting
$L:=L_1+L_2+p$ finishes the proof.
\end{proof}

\begin{rmk}
Thus, even if the normalization of a word is not finite, you can
add a ``period killer" of a fixed length so that it will become
such. This is in fact the only property we will be using. It
looks weaker than WF and we have been even tempted to call it
something like PWF (positively weakly finitary); the reason for
not doing so is the fact that there are no examples of Pisot
numbers that are PWF but not WF (actually, as we already
mentioned above, there are no examples of Pisot numbers that are
not WF at all!).
\end{rmk}

\begin{lemma}
For $p_d^+$-a.e. sequence $x\in\Si_d^+$ there exists $n$ such that
$\n_{\be,d}^+(x_1\dots x_n)$ is finite. \label{blk}
\end{lemma}
\begin{proof} The proof is similar to the one of
\cite[Proposition~18]{S2}. Let
\begin{equation}
\mathfrak{A}_n=\left\{x\in\Si_d^+ : \n_{\be,d}^+(x_1\dots
x_n)\,\,\mathrm{is\ finite}\right\}\label{An}
\end{equation}
and $\mathfrak{B}_n=\Si_d^+\setminus\mathfrak{A}_n$. Our goal is
to show that there exists a constant $\ga=\ga(\be,d)\in(0,1)$ such
that
\begin{equation}
p_d^+\left(\bigcap_{k=1}^n\mathfrak B_k\right)\le\ga^n,\quad
n\ge1. \label{ineq}
\end{equation}
Let $L=L(\be,d)$ be as in Lemma~\ref{L}. We have
\begin{eqnarray*}
p_d^+\left(\bigcap_{k=1}^n\mathfrak B_k\right)  &\le
&\prod_{k=2}^n p_{d}^+(\mathfrak B_k\mid \mathfrak B_{k-1}\cap
\dots \cap \mathfrak
B_1) \\
&\le &\prod_{k=2}^{[n/L]}p_d^+(\mathfrak B_{Lk}\mid \cap
_{j=1}^{Lk-L-1}\mathfrak B_j).
\end{eqnarray*}
By Lemma~\ref{L},
\[
p_d^+(\mathfrak A_{Lk}\mid\mathfrak E )\ge d^{-L},
\]
for any $\mathfrak E$ in the sigma-algebra generated by
$x_{Lk-L-1},\dots,x_1$. Hence
\[
p_d^+(\mathfrak B_{Lk}\mid\cap_{j=1}^{Lk-L-1}\mathfrak B_j)\le
1-d^{-L},
\]
and
\[
p_d^+\left(\bigcap_{k=1}^n\mathfrak
B_k\right)\le(1-d^{-L})^{[n/L]}.
\]
It suffices to put $\gamma:=(1-d^{-L})^{1/2L}$, which proves
(\ref{ineq}) and the lemma.
\end{proof}

\begin{cor}\label{maincor} For $p_d^+$-a.e. sequence $x\in\Si_d^+$
its normalization is blockwise, i.e., there exists a sequence
$(n_k)_{k=0}^\infty$ such that $n_0=0$ and
\[
\n^+_{\be,d}(x)=\n^+_{\be,d}(x_1\dots
x_{n_1})\n^+_{\be,d}(x_{n_1+1}\dots x_{n_2})\dots
\]
(a concatenation of finite words), and $\n_{\be,d}(x_{n_k+1}\dots
x_{n_{k+1}})$ is finite of the length $n_{k+1}-n_k$ for any
$k\ge0$.
\end{cor}
\begin{proof}
Recall that by the well-known result from \cite{FrSo} quoted in
Lemma~\ref{L}, there exists a number $K=L_1\in\BN$ such that if
the normalization of a word of length~$n$ in a Pisot base with a
fixed alphabet is finite, then the length of its normalization is
at most $n+K$ (here $K$ depends on $\be$ and the alphabet only).
This means that if $w_1, w_2$ are two words with finite
normalizations, then the normalization of $w_10^{2K}w_2$ is the
concatenation of the normalizations of $w_10^K$ and $0^Kw_2$.

Set
$$
\mathfrak D=\{x\in\Si_d^+\mid \exists (l_k)_1^\infty: x_j\equiv0,\
l_k\le j\le l_k+2K, \forall k\ge1 \}.
$$
Obviously, $p_d^+(\mathfrak D)=1$. Put
$$
\mathfrak A^+=\mathfrak D\cap\bigcup_{n=1}^\infty \mathfrak A_n,
$$
where $\mathfrak A_n$ is given by (\ref{An}). Still we have
$p_d^+(\mathfrak A^+)=1$. Thus, the probability that at the end of
the first block in Lemma~\ref{blk} there are $2K$ consecutive
zeros, is also 1. Therefore, by the above, the normalizations of
the first block and all the rest will be totally independent.
Consider the second block, then the third one, etc. - and then
take the countable intersection of all the sets obtained. This is
a sought set of full measure $p_d^+$.
\end{proof}

\medskip\noindent\textbf{3. Two-sided Erd\"os measure and
conclusion of the proof}. Now we are ready to define the
two-sided normalization and -- consequently -- the two-sided
Erd\"os measure.

\begin{Def} Let $x=(x_n)_{-\infty}^\infty\in\Si_d$. Put
\begin{equation}
\n_{\be,d}(x):=\lim_{N\to+\infty}\n^+_{\be,d}(x_{-N},x_{-N+1},
\dots), \label{norm}
\end{equation}
where the limit is taken in the natural (weak) topology of
$\Si_d$. By the previous corollary and the fact that the measure
$p_d$ is the weak limit of the measures $S_d^n(p_d^+)$ (where
$S_d$ denotes the shift on $\Si_d$), we conclude that the map
$\n_{\be,d}:\Si_d\to X_\be$ is well defined and blockwise (in the
sense of the previous corollary) for $p_d$-a.e. sequence
$x\in\Si_d$. We will call it the {\em two-sided normalization} (in
base $(\be,d)$).
\end{Def}

\begin{rmk}
It is worth noting that if $\be$ is an algebraic unit (i.e., if
$\be^{-1}\in\BZ[\be]$), then there exists an alternative way of
defining $\n_{\be,d}$ via the torus. Namely, let
$x^m=k_1x^{m-1}+\dots+k_m$ be the characteristic equation for
$\be$ ($k_m=\pm1$) and $\BT^m=\BR^m/\BZ^m$. Let $T_\be$ be the
automorphism of $\BT^m$ determined by the companion matrix
$M_\be$ for $\be$, i.e.,
\[
M_{\beta}=\left(
\begin{array}
[c]{ccccc}%
k_{1} & k_{2} & \ldots &  k_{m-1} & k_{m}\\
1 & 0 & \ldots & 0 & 0\\
0 & 1 & \ldots & 0 & 0\\
\ldots & \ldots & \ldots & \ldots & \ldots\\
0 & 0 & \ldots & 1 & 0
\end{array}
\right).
\]
Let $\mathcal{H}(T_\be)$ denote the group of points homoclinic
to zero, i.e., $\mathbf{t}\in\mathcal{H}(T_\be)$ iff
$T_\be^n(\mathbf{t})\to0$ as $n\to\pm\infty$. A homoclinic point
$\mathbf{t}$ is called \textit{fundamental} if the linear span
of its $T_\be$-orbit is the whole group $\mathcal{H}(T)$. It is
well known that such points always exist for $T_\be$ (actually
they exist for any automorphism of $\BT^m$ which is
$SL(m,\BZ)$-conjugate to $T_\be$ -- see \cite{Ver, S2}). Now let
the map from $X_\be$ onto $\BT^m$ be defined as follows:
\begin{equation}
F_{\mathbf t}(\e)=\sum_{n\in\BZ}\e_nT_\be^{-n}(\mathbf t),
\label{F}
\end{equation}
where $\mathbf t$ is fundamental. It is easy to show that the
series does converge on the torus \cite{Sch2000, S2} whenever
the $\e_n$ are bounded. Let $\tau_\be$ denote the shift on
$X_\be$, i.e., $\tau_\be(\e)_n=\e_{n+1}$. In \cite{S2} it is
shown that if $\be$ is WF, then $F_{\mathbf t}$ is one-to-one
a.e. and conjugates the shift $\tau_\be$ and $T_\be$.

Now we state without proof that similarly to the one-sided
normalization (see (\ref{n+})), the two-sided normalization can be
computed by the formula
\[
\n_{\be,d}=F_{\mathbf{t}}^{-1}\wt F_{\mathbf{t},d},
\]
where the projection $\wt F_{\mathbf{t},d}:\Si_d\to\BT^m$ is
given by the same formula~(\ref{F}) as $F_{\mathbf{t}}$ with
$\e_n$ replaced by $x_n$. In particular, it is well defined a.e.
and does not depend on a choice of $\mathbf t$.
\end{rmk}

\begin{Def} The projection $\nu_{\be,d}:=\n_{\be,d}(p_d)$ is
called the {\em two-sided Erd\"os measure}.
\end{Def}

We have the following diagram:

\[
\begin{CD}
\Si_d @>{S_d}>> \Si_d \\
@V{\n_{\be,d}}VV @VV{\n_{\be,d}}V \\
X_\be @>{\tau_\be}>> X_\be
\end{CD}
\]

Since the the two-sided normalization obviously commutes with the
shift (see (\ref{norm})), the diagram commutes as well, whence
the two-sided Erd\"os measure is also shift-invariant (unlike the
one-sided Erd\"os measure!). This is because
\[
\tau_\be(\nu_{\be,d})=\tau_\be\n_{\be,d}(p_d)=\n_{\be,d}S_d(p_d)=
\n_{\be,d}(p_d)=\nu_{\be,d}.
\]
Moreover, since the automorphism $(\Si_d,p_d,S_d)$ is Bernoulli,
so is the automorphism $(X_\be,\nu_{\be,d},\tau_\be)$ -- by
Ornstein's Theorem which claims that all Bernoulli factors are
Bernoulli \cite{Orn}.

Let $\rho_d:\Si_d\to\Si_d^+$ and $\rho_\be:X_\be\to X_\be^+$
denote the natural projections. Then
\[
\mu_{\be,d}^+ = (\n_{\be,d}^+\rho_d)(p_d).
\]
Let
\[
\nu_{\be,d}^+ :=
(\rho_\be\n_{\be,d})(p_d)=(\phi_\be\rho_\be)(\nu_{\be,d}).
\]
By the above, $\nu_{\be,d}^+$ is $\tau_\be^+$-invariant and its
natural extension is Bernoulli.

Note that one of the reasons why one may have difficulties with
the one-sided Erd\"os measure is because the operations of
normalization and projection do not commute (and therefore,
$\mu_{\be,d}^+$ is not shift-invariant). However, in a sense
these operations are ``commuting up to a finite number of
coordinates", which allows us to finish the proof of
Theorem~\ref{main}.

\begin{lemma} The measures $\mu_{\be,d}^+$ and $\nu_{\be,d}^+$ are
equivalent.
\end{lemma}
\begin{proof} Let $P_{\be,d}=\rho_\beta\circ\n_{\be,d}$
and $Q_{\be,d}=\n_{\be,d}^+\circ\rho_d$. Since $p_d$ is
preserved by the action of the group that changes a finite
number of coordinates, it suffices to show that there exist two
maps $C:\Si_d\to\Si_d$ and $C':\Si_d\to\Si_d$ defined
$p_d$-almost everywhere with the following properties: each of
them is a step function with a countable number of steps, it
changes just a finite number of coordinates of $x$ and also
\begin{equation}
P_{\be,d}(C(x))=Q_{\be,d}(x),\,\,
Q_{\be,d}(C'(x))=P_{\be,d}(x).\label{change}
\end{equation}
If we construct such functions, this will prove
Theorem~\ref{main}, because then we will have
\[
(\n_{\be,d}^+\rho_d)(p_d)\prec(\rho_\be\n_{\be,d})(p_d),\,\,
(\rho_\be\n_{\be,d})(p_d)\prec(\n_{\be,d}^+\rho_d)(p_d),
\]
i.e., $\nu_{\be,d}^+\approx\mu_{\be,d}^+$.

Let $\mathfrak A$ be the two-sided analog of $\mathfrak A^+$
defined in the proof of Corollary~\ref{maincor}, namely,
$\mathfrak A$ is the set of all sequences in alphabet $\A_d$
whose normalization $\n_{\be,d}$ is blockwise in the sense of
Lemma~\ref{blk}. This set has full measure $p_d$. Let
$x\in\mathfrak A$; then $x$ can be represented in the block form
$x=(\dots B_{-2}B_{-1}B_0B_1\dots)$ and
\[
\n_{\be,d}(x)=(\dots
\n_{\be,d}(A_{-1})\n_{\be,d}(A_0)\n_{\be,d}(A_1)\dots),
\]
where each word $\n_{\be,d}(A_n)$ is of the same length as
$A_n$.

Let $A_0=(x_{-a}\dots x_{b})$ with $a>0, b>0$ (one can always
achieve this by merging blocks). By the above, we have
$P_{\be,d}(x)=(\e_1,\dots,\e_b,*)$ and
$Q_{\be,d}(x)=(\e'_1,\dots,\e'_b,*)$ (where the star indicates
one and the same tail), i.e., the difference is only at the
first $b$ places. Thus, it is easy to guess what $C$ and $C'$
may look like. Namely, put
\[
(C(x))_j=
\begin{cases} x_j,&j<-a\,\,\mathrm{or}\,\, j>b\\
0,& -a\le j\le 0\\
\e_j',&1\le j\le b
\end{cases}
\]
and
\[
(C'(x))_j=
\begin{cases} x_j,&j<-a\,\,\mathrm{or}\,\, j>b\\
0,& -a\le j\le 0\\
\e_j,&1\le j\le b.
\end{cases}
\]
Both functions are obviously well defined for $p_d$-a.e. $x$,
are step functions with a countable number of steps and change a
finite number of coordinates. The equalities in (\ref{change})
are satisfied as well, which proves Theorem~\ref{main}.
\end{proof}

\medskip\noindent\textbf{4. Acknowledgement.} The author wishes to
thank E.~Olivier and A.~Thomas for helpful discussions and
suggestions.

\end{document}